\documentclass[english, 11pt]{article}
\usepackage[]{fontenc}
\usepackage{fancyhdr}
\usepackage[nottoc]{tocbibind}
\usepackage{babel}
\usepackage{amsmath}
\usepackage{amscd}
\usepackage{amsthm}
\usepackage{amsfonts}
\newtheorem{lem}{Lemma}[section]
\newtheorem{thm}[lem]{Theorem}
\newtheorem{pro}[lem]{Proposition}
\newtheorem{cor}[lem]{Corollary}

\newtheorem*{theo}{Theorem}
\theoremstyle{definition}
\newtheorem{rem}[lem]{Remark}
\theoremstyle{definition}
\newtheorem{df}[lem]{Definition}
\theoremstyle{remark}
\newtheorem*{cl}{Claim}
\newtheorem{st}{Step}

\numberwithin{equation}{section}

\newcommand{\M}{\mathcal{M}}

\newcommand{\Ml}{M^{loc}}
\newcommand{\Mls}{M_s ^{loc}} 
\newcommand{\Ms}{M_s}
\newcommand{\Mle}{M_{\eta} ^{loc}} 
\newcommand{\OO}{\mathcal{O}}
\newcommand{\Spec}{\text{Spec}}
\newcommand{\frob}{\text{Frob}}
\newcommand{\Tr}{\text{Tr}}
\newcommand{\la}{\langle}
\newcommand{\ra}{\rangle}

\newcommand{\rk}{\text{rk}}
\newcommand{\spa}{\text{sp}}

\begin{document}

\title{Local models for ramified unitary groups}
\author{Nicole Kr\"amer\\TU Berlin - Department of Computer Science\\Franklinstr. 28/29\\D-10587 Berlin\\nkraemer@cs.tu-berlin.de}
\maketitle

\section{Introduction}
Given a Shimura variety $S=S(G,h)$ over a number field $E$  it is useful to have a model over the ring of integers $\mathcal{O}_E$\,, i.e.\ a scheme $S' \rightarrow \Spec (\mathcal{O}_E)$ such that the diagram
\begin{eqnarray*}
\begin{CD}
S @>>> \Spec(E) \\
@VVV   @VVV \\
S' @>>>\Spec (\mathcal{O}_E)  
\end{CD}
\end{eqnarray*}
is cartesian. To qualify as a model, the scheme $S'$ should  at least be flat over $\Spec (\mathcal{O}_E)$ and have only ``good'' singularities. In the PEL-case, $S$ can be written as a moduli space of abelian varieties with polarization, an action by a semisimple $\mathbb{Q}$-algebra and a $K$-level strucure where $K \subset G(\mathbb{A}_f)$ is an open and compact subgroup. If we fix a prime $p$\,, the investigation of a model $S'$ at $p$ can be reduced to the investigation of the associated ``local model'' $M$, introduced in \cite{rz} -- at least if $K=K^p.K_p$ where $K^p \subset G(\mathbb{A}^p _f)$ and $K_p \subset G(\mathbb{Q}_p)$ is a parahoric subgroup.\\
In this work we will consider the case that $K_p$ is a maximal parahoric subgroup of $G(\mathbb{Q}_p)$\,. We further fix a prime $\mathfrak{p}$ of $\mathcal{O}_E$ lying over $p$ and denote by $R$ the completion of $(\mathcal{O}_E)_{\mathfrak{p}}$\,. Then the local model $M$ is a closed subscheme of a Grassmannian, defined over   $\Spec (R)$\,. Locally for the \'etale topology $M$ coincides with $S'\otimes_{\mathcal{O}_E} R$\,. Hence $S'$  is flat if and only if the local model is and both schemes have the same singularities.\\
We are going to study  local models $M_{r,s}$ of Shimura varieties where  $G$ is a reductive group over $\mathbb{Q}$  such that $G \otimes_{\mathbb{Q}} \mathbb{R}=GU(r,s)$\,. In \cite{rz} it is conjectured that $M_{r,s}$ is flat over the base $R$. This is true if $p$ is unramified by the work of G\"ortz \cite{G2}, but an argument of Pappas shows that this is not true in general: Using the facts that the generic fiber of $M$ is a Grassmannian of dimension $rs$ (cf. \cite{Pappas} or proposition \ref{generic fiber} below) and that the special fiber does not depend on the partition $(r,s)$ of $n=r+s$ we show nonflatness in the case $|r-s| > 1$ (see corollary \ref{notflat} for a somewhat different proof of Pappas' statement).  Hence these ``naive'' local models need to be modified. Therefore we investigate  a closed subscheme $M_{r,s} ^{loc}$ of the naive local model -- introduced by Pappas -- in section \ref{local model}, which stands a better chance to be flat. In fact, Pappas is able to show the flatness of $M_{r,s} ^{loc}$ if $(r,s)=(n-1,1)$ (see \cite{Pappas}, theorem 4.5). In the case  $|r-s|\leq 1$\,, Pappas conjectures that $M_{r,s}=M_{r,s} ^{loc}$ and is able to prove this in the case $(r,s)=(2,1)$.\\
Another interesting task in the theory of Shimura varieties  is to calculate the local factors of the semisimple Hasse-Weil zeta function at places where the Shimura variety has bad reduction. (For a definition of the semisimple zeta function see \cite{HN}). It is related to the alternating semisimple trace  
\begin{align*}
\text{Tr}^{\text{ss}}(\text{Frob}_q;(R\Psi \mathbb{Q}_l)_x ^I) 
\end{align*} 
of the geometric Frobenius on the invariants of the inertia group of the sheaf of nearby cycles.  We are going to calculate this trace in the case $(r,s)=(n-1,1)$\,. Again we can restrict to the local model associated to the Shimura variety. We cannot calculate the alternating semisimple trace  directly, as the local model does not have  semi-stable reduction.  We have to proceed in several steps: First we have to resolve the  singularities of the local model. In \cite{Pappas}, this was done by blowing up the singular locus. Our approach is different:  We define a  resolution $\tau: \mathcal{M} \rightarrow M_{n-1,1}^{loc}$ by posing a moduli problem analogous to  the Demazure-resolution of a Schubert variety in the Grassmannian. It was suggested by Pappas and Rapoport. The following result is the main result of this article. It was conjectured by Pappas and Rapoport.
\begin{theo}
Let $(r,s)=(n-1,1)$\,. $\M$ is a regular scheme of dimension $n$\,.Its special fiber is reduced and has normal crossings. The  morphism $\tau:\M \rightarrow M^{loc}$ is proper and induces an isomorphism outside a point $y \in M^{loc}$. The fiber of $\tau$ over $y$ is a projective space of dimension $n-1$\,. The generic fiber of $\M\rightarrow \Spec(R)$ is smooth and geometrically irreducible of dimension $n-1$\,. The special fiber is  the union of two smooth irreducible varieties of dimension $n-1$ -- one of them being the fiber over $y$ -- which intersect in a smooth irreducible variety of dimension $n-2$\,. The morphism $\M \rightarrow \Spec(R)$ is flat.
\end{theo}
As  $\mathcal{M}$ has semistable reduction, we are  able to compute the alternating semisimple trace on $\mathcal{M}$ which then can be related to the alternating semisimple trace on $M_{n-1,1} ^{loc}$\,.\\
\\
In conclusion, I would like to express my gratitude to those who helped me with this work. First of all, I would like to thank M. Rapoport for introducing me into this area of mathematics and his steady interest in my work. I am very grateful to T. Wedhorn. His mathematical advice, his encouragement and his infinite patience accompanied me during my work on this article. Finally, I would like to thank U. G\"ortz for helpful remarks on an earlier version of this work.
\section{The naive local model}
\label{naive}
Let $R_0$ be a complete discrete valuation ring with uniformizing element $\pi_0$ and perfect residue field $k$\,. We assume  that $\text{char}\, k \not=2$\,. Set $K_0 =\text{Quot}(R_0)$ and
\begin{align*}
K=K_0[X]/(X^2-\pi_0)\;.
\end{align*}
Denote by $\pm \pi$ the roots of $X^2-\pi_0$ in $K$\,. The ring $R= \OO_K$ of integers is a discrete valuation ring with uniformizing element $\pi$ and residue field $k$\,. For $n\geq 3$ set $V=K^n$ with canonical basis $e_1,\ldots,e_n$ and denote by $\Lambda$ the standard lattice $\bigoplus R e_i$ in $V$\,. We write $\overline{\Lambda}$ for $\Lambda \otimes k$\,. The uniformizing element $\pi$ induces a $R_0$ - linear mapping on $\Lambda$ which we denote by $\Pi$\,. 
Define on $V$ an alternating $K_0$-bilinear mapping
\begin{align*}
\la \; ,\, \ra \colon &V \times V \rightarrow K_0
\end{align*}
by imposing that for $i,j \in \{1,\ldots,n\}$
\begin{align*}
\la e_i,e_j\ra =0\,, \la e_i,\Pi e_j \ra =\delta_{ij}\;.
\end{align*}
For any element $a=x + \pi y \in R$ $(x,y \in R_0)$ set $\overline{a}=x - \pi y$\,. Then we have
\begin{align*}
\la av,w \ra=\la v,\overline{a}w\ra 
\end{align*}
for $a \in R$ and $v,w \in V$\,. Fix natural numbers $r,s$ with $r+s=n$\,.
\begin{df}
Let $M_{r,s}$ be the functor  which associates to each scheme $S$ over $\Spec (R)$ the set of subsheaves $F$ of $(R \otimes_{R_0} \OO_S)$-modules of $\Lambda \otimes_{R_0} \OO_S$ such that
\begin{enumerate}
\item $F$ as an $\mathcal{O}_S$-module is locally on $S$ a direct summand of rank $n$\;;
\item $F$ is totally isotropic for $\la \;,\;\ra \otimes_{R_0} \OO_S$\;;
\item \label{char}$\text{char}_{\Pi|F}(T)=(T+ \pi)^r (T- \pi)^s $\;.
\end{enumerate}
\end{df}
The functor $M_{r,s}$ is represented by a closed subscheme of $Gr(n,2n) \otimes R$ hence $M_{r,s}$ is a projective $\Spec (R)$-scheme. (Here we denoty by $Gr(a,b)$ the scheme parameterizing locally direct summands of rank $a$ of a free module of rank $b$).\\
The next result is due to Pappas (cf. \cite{Pappas}, proposition 3.8). Recall that  $K=\text{Quot}(R)$\,.
\begin{pro}[Pappas]
\label{generic fiber}
We have 
\begin{align*}
M_{r,s} \times_{\Spec(R)} \Spec(K) \simeq Gr(n-s,n) \otimes K \;.
\end{align*}
In particular, the generic fiber of $M_{r,s}$ is  smooth and geometrically irreducible of dimension $r s$\,.
\end{pro}
\begin{proof}
This proof is more or less a replication of Pappas' proof in \cite{Pappas}. Let $B$ be a $K$-algebra. We have an isomorphism of $B$-algebras
\begin{align*}
R\otimes_{R_0} B &\rightarrow B \oplus B\;, \\
a \otimes b &\mapsto (ab,\overline{a}b)\;.
\end{align*}
Denote by $f_0,f_1$ the two idempotent elements of $R
\otimes_{R_0} B$ corresponding to $(1,0)$ and $(0,1)$ in  $ B
\oplus B$\,. Set 
\begin{align*}
G_i=f_i (R \otimes_{\mathbb{Z}_p} B)^n\,, i=0,1\;.
\end{align*}
Then $G_i$ is totally isotropic with respect to $\la\;, \;\ra \otimes B$\,. 
For $i \in \mathbb{Z}/2\mathbb{Z}$\,, consider the map
\begin{align*}
\psi_i \colon & G_i \longrightarrow  Hom_B(G_{i+1},B)\;,\\
       & x  \longmapsto  \la x, \;\ra |_{G_{i+1}}\;.
\end{align*}
As $\la \;, \;\ra$ is perfect and $G_i$ is isotropic, $\psi _i$ is an
isomorphism. Now consider the functor $N$ which associates to each
$K$-algebra $B$ the set of pairs $(F_0,F_1)$ of $B$-modules which are a direct summand of $B^n$ of rank $s$ (resp. $r$) such that $\psi_i
(F_i)$ vanishes on $F_{i+1}$ ($i \in \mathbb{Z}/2\mathbb{Z}$). (Here we view $F_i$ as subsets of $G_i \simeq B^n$).
We get a bijection
\begin{align*}
M_{r,s}(B) \rightarrow N(B)
\end{align*}
via
\begin{align*}
F\mapsto (f_0F,f_1F)\;.
\end{align*}

On the other hand, the pair $(F_0,F_1)$ is determined by $F_0$\,. Indeed,  we
 have $F_1 \subseteq F_0 ^{\bot}$\,.  This is an inclusion of direct summands of $B^n$ of rank $r=n-s$\,, hence they must be equal. We conclude that
\begin{align*}
M_{r,s}(B)= Gr(n-s,n)(B)\,.
\end{align*} 
\end{proof}
\begin{cor}[Pappas]
\label{notflat}
We have
\begin{align*} 
\dim M_{r,s}\otimes_R k \geq
\begin{cases} 
\frac{n^2}{4} & \text{if }n \text{ is even}\\ 
\frac{n^2 -1}{4} & \text{if }n \text{ is odd\;.}
\end{cases}
\end{align*}
In particular, $M_{r,s}$ is not flat if $|r-s| > 1$\,.
\end{cor}
\begin{proof}
 As $M_{r,s} \rightarrow \Spec(R)$ is proper  we have the inequality  
\begin{align*}
\dim M_{r,s}\otimes_R k \geq \dim M_{r,s}\otimes_R K\;.
\end{align*}
On the other hand, the special fiber of $M_{r,s}$ does not depend on the partition $(r,s)$ of $n$\,. Indeed, condition \ref{char} is just $\text{char}_{\Pi|F}=T^n$ in the special fiber.
Using this and proposition \ref{generic fiber} we see that  $\dim M_{r,s}\otimes_R k \geq r's'$ for every pair $(r',s')$ with $r'+s'=n$\,. Choosing $(r',s')$ such that $r's'$ is maximal gives the inequality. In particular, if $|r-s| > 1$\,, the dimension of the two fibers are not equal, hence the proper morphism $M_{r,s} \rightarrow \Spec(R)$ cannot be flat.
\end{proof}
\section{The local model}
\label{local model}
\begin{df}
\label{df:lm}
If $r=s$\,, set $M^{loc} _{r,s}= M_{r,s}$\,. If $r \neq s$ let $M^{loc} _{r,s}$ be the functor $(Sch/\Spec R) \longrightarrow (Sets)$ which associates to each scheme $S$ over $\Spec(R)$ the subset of $M_{r,s}(S)$ which consists of subsheaves $F$ such that
\begin{enumerate}
\item \label{(a)} $\bigwedge^{r+1}(\Pi -\pi|F)=(0)$\;;
\item \label{(b)}$\bigwedge^{s+1}(\Pi +\pi|F)=(0)$\;.
\end{enumerate}
\end{df}
The functor $M^{loc} _{r,s}$ is represented by a closed subscheme of $M_{r,s}$\,. We fix $(r,s)$ and set $M=M_{r,s}$ and $\Ml = M^{loc} _{r,s}$\,. We denote by $\Ml_{\eta}$ and $\Ml_s$ the generic and special fiber of $\Ml$ respectively.\\
While the naive local model $M$ cannot be flat over $\Spec (R)$ if $|r-s|>1$\,, Pappas conjectures  in \cite{Pappas} that $\Ml$ is flat for every pair $(r,s)$\,. He is able to prove the following. Denote by $y$ the $\Spec (k)$-valued point $\Pi \overline{\Lambda}$ of $\Mls$\,.
\begin{thm}[Pappas]
\label{pappas} Suppose that $(r,s)=(n-1,1)$\,.  Then $\Ml$ is normal and Cohen-Macaulay and the morphism $\Ml \to \Spec (R)$ is flat and and it is smooth outside the closed subscheme defined by $y$\,.
\end{thm}
\begin{proof}
See \cite{Pappas}, theorem 4.5 .
\end{proof}
\begin{rem}
We need only one consequence of theorem \ref{pappas} in this work: We  use the fact  that $\Ml_{n-1,1}$ is smooth outside the point $y$ in proposition \ref{red}. 
\end{rem}
\begin{pro}
The generic fibers of $\Ml$ and $M$ coincide, in particular the generic fiber of $M^{loc}$ is a smooth, geometrically irreducible variety of dimension $rs$\,.\end{pro}
\begin{proof}
In the proof of proposition \ref{generic fiber} we saw  that for every $K$-algebra $B$ and every
$F \in M(B)$ there is a decomposition $F=F_0 \oplus F_1$ with $\rk F_0=s$ and  $\rk F_1=r$  such that $\Pi$
acts on $F_0$ via $\pi$ and on $F_1$ via $-\pi$\,. Hence
\begin{align*}
\sideset{}{^{r+1}} \bigwedge (\Pi -\pi|F)&=\sideset{}{^{r+1}} \bigwedge (+2 \pi|F_1) =(0)\;,\\
\sideset{}{^{s+1}} \bigwedge (\Pi +\pi|F)&=\sideset{}{^{s+1}} \bigwedge (-2 \pi|F_0) =(0)\;.
\end{align*}
\end{proof}
\begin{rem}
We saw in corollary \ref{notflat} that $\Ml\not=M$ if  $|r-s|> 1$. Pappas conjectures in \cite{Pappas} that $\Ml=M$ if $|r-s|\leq 1$ and is able to prove this in the case $(r,s)=(2,1)$ (see \cite{Pappas} 4.15). 
\end{rem}

\section{A resolution of the singularities}\label{resolution}
In this section we introduce  a resolution of the singularities of the local model in the case $(r,s)=(n-1,1)$\,.
\begin{df}
Let $\M$ be the functor $(Sch/\Spec (R)) \longrightarrow (Sets)$ which associates to each scheme $S$ over $\Spec (R)$ the set of pairs $(F_0,F)$ where $F_0$ and  $F$ are $\OO_S$-submodules of $\Lambda {\otimes}_{R_0} \OO_S$ such that 
\begin{enumerate}
\item \label{eins}$F$ as an $\mathcal{O}_S$-module is locally on $S$ a direct summand of rank $n$\;; 
\item \label{zwei}$F_0$ as an $\mathcal{O}_S$-module is locally on $S$ a direct summand of rank $1$\;;
\item \label{drei} $F_0 \subseteq F$\;;
\item \label{vier}$F$ is totally isotropic for $\la\;,\;\ra {\otimes}_{R_0} \OO_S$\;;
\item \label{fuenf}$(\Pi + \pi)F \subseteq F_0$\;;
\item \label{sechs}$(\Pi - \pi)F_0 =(0)$\;.
\end{enumerate}
\end{df} 
The functor is represented by a projective $\Spec (R)$-scheme. We denote by $\M_{\eta}$ and $\M_s$ the generic and special fiber of $\M$ respectively.\\
Consider the  projection morphism 
\begin{align*}
\tau \colon & \M \longrightarrow \Ml,
\end{align*}
which is given by 
\begin{align*}
(F_0,F)\mapsto F 
\end{align*}
on $S$-valued points.
The map $\tau$ is well defined: Condition \ref{drei} and \ref{fuenf} imply that $F$ is $\Pi$-stable. If we assume that $F$ is a free module
spanned by elements $v_1,\ldots, v_n$ with $F_0=\spa (v_1)$ the map $\Pi$ is
represented by the matrix
\begin{align*}
\begin{pmatrix}
\pi&\lambda_2  &\hdots&\lambda_n \\
0 &-\pi&\hdots&0 \\
 \vdots&        &        &  \\
0&0     &  \hdots&-\pi
\end{pmatrix}
\end{align*}
with suitable $\lambda_2,\ldots,\lambda_n$\,. We conclude that the condition on $\text{char}_{\Pi|F}$ is satisfied, that $\bigwedge^{n}(\Pi - \pi|F)=(0)$ and that $\bigwedge^{2}(\Pi +\pi|F)=(0)$\,.
As $\M$ and $\Ml$ are proper $\Spec(R)$-schemes, the morphism $\tau$ is proper.
\begin{rem}
\label{rem:generic}
The morphism  $\tau$ induces an isomorphism on the generic fibers. Indeed, the proof of proposition \ref{generic fiber} shows that for every $K$-algebra $B$ and every $F \in \Ml(B)$ we have a decomposition $F=F_0 \oplus F_1$ and $F\mapsto (F_0,F)$ defines the inverse map to $\tau$\,.
\end{rem}
Recall that we denote  by $y$ the $\Spec (k)$-rational point $\Pi \overline{\Lambda}$\,.
\begin{pro}\label{red}
 The map $\tau$ induces an isomorphism
\begin{align*}
\tau' \colon \M \setminus \tau^{-1}(y) \to \Ml \setminus \{y\}\;.
\end{align*}
\end{pro}
\begin{proof}
We know from remark \ref{rem:generic} that $\tau'$ is an isomorphism on the generic fibers. Next consider $\tau'$ on the special fibers. A point $z \in \Ml_s \setminus \{y\}$ corresponds to a $n$-dimensional $\kappa(z)$-subvectorspace $F$ of $\Lambda_{\kappa(z)}= \Lambda \otimes \kappa(z)$  such that $F \not= \Pi \Lambda_{\kappa(z)}$\,. It follows that $\Pi F \not= (0)$\,. Let $T$ be a $\Spec (\kappa(z))$-scheme and  $(F_0, F\otimes_{\kappa(z)} T)$ a $T$-valued point of the fiber in $z$\,. Set $F_T=F\otimes_{\kappa(z)} T$\,. By definition we have $\Pi F_T \subseteq F_0$\,. Take a point $t \in T$\,. As the morphism $\Spec(\kappa(t)) \rightarrow \Spec (\kappa(z))$ is faithfully flat we have $(\Pi F)_{\kappa(t)} \not=0$\,, hence it must equal $(F_0)_{\kappa(t)}$\,. By Nakayama's Lemma $(F_0)_t=(\Pi F_T)_t$\,, hence  $F_0=\Pi F_T$\,. We conclude that the scheme-theoretic fiber in any point $z$ of $\Ml \setminus \{y\}$ equals $\Spec(\kappa(z))$\,. It follows  that $\tau'$  is an isomorphism: As $\tau$ is proper, the morphism $\tau'$ is also proper. Furthermore we have just seen that it is quasifinite, hence Zariski's main theorem implies that $\tau'$ is finite. As $\Ml\setminus\{y\}$ is reduced (cf. proposition \ref{pappas}), \cite{bourbaki} \S3 $n^02$ proposition 7  implies that $\tau'$ is flat.  Lemma \ref{pro:gleichheit} gives the result.
\end{proof}
\begin{lem}

\label{pro:gleichheit}
Let $S$ be a scheme and $f \colon X \rightarrow Y$ be a finite and flat
morphism of $S$-schemes. Assume that $f$ is an isomorphism on the fibers. Then $f$ is an isomorphism.
\end{lem}
\begin{proof}. It suffices to show the statement in the case $X=\Spec (B)$ and  $Y=\Spec (A)$\,. As $B$ is a finitely generated, flat $A$-module, $B$ is locally free. We may assume that $B$ is free, say $B \simeq A^n$\,. As the morphism is an isomorphism on the fibers, we conclude that for a prime ideal $p \in A$ we have $B \otimes_A \kappa (p) \simeq \kappa (p)^n \simeq \kappa (p)$\,, hence $n=1$\,.  
\end{proof}
\begin{thm}
\label{thm:reg}
$\M$ is a regular scheme of dimension $n$ with reduced special fiber. The generic fiber of $\M$ is smooth and geometrically irreducible of dimension $n-1$\,. The special fiber is  the union of two smooth irreducible varieties of dimension $n-1$ which intersect in a smooth irreducible variety of dimension $n-2$\,. The morphism $\M \rightarrow \Spec(R)$ is flat.
\end{thm}
\begin{rem}
In the proof of theorem \ref{thm:reg} we will see that the second part of the  theorem  is also correct for the geometric special fiber of $\M$\,.
\end{rem}
\begin{proof}[Proof of theorem \ref{thm:reg}]
It follows from remark \ref{rem:generic} that the generic fiber of $\M$ is smooth and geometrically irreducible of dimension $n-1$ as $\Mle$ has these properties. Now consider the special fiber of $\M$\,. We proceed in several steps.
\begin{st}
Set
\begin{align*}
Z_1 =\tau_s ^{-1}(y) = y \times_{\Ml_s} \M_s\;.
\end{align*}
We have $Z_1 \simeq \mathbb{P}(\Pi \overline{\Lambda})$\,.
Indeed, for every $\Spec(k)$-scheme $S$ and for every  $F_0 \in \mathbb{P}(\Pi \overline{\Lambda})(S)$ the pair  $(F_0,\Pi \overline{\Lambda}_S)$ satisfies the conditions \ref{fuenf} and \ref{sechs}.
\end{st}
\begin{st}
Next, consider the morphism
\begin{align*}
\varphi \colon & \mathcal{M}_s \longrightarrow \mathbb{P}(\Pi \overline{\Lambda})\;,
\end{align*}
which is given by
\begin{align*}
(F_0,F)\mapsto F_0
\end{align*}
on $S$-valued points ($S$ a $\Spec (k)$-scheme).
Define a  form  on $\Pi \overline{\Lambda} $ by
\begin{align*}
\{\Pi v,\Pi w\}=\la \Pi v, w\ra \;.
\end{align*}
This is well-defined: If $\Pi w=\Pi w'$ then $w=w' + \Pi u$ and $\{\Pi v,\Pi w\}=\la \Pi v, w\ra =\la \Pi v,w'\ra + \la \Pi v,\Pi u \ra=\{\Pi v,\Pi w\}$\,. The form is symmetric as
\begin{align*}
\{\Pi v,\Pi w\}=\la \Pi v, w\ra=-\la w,\Pi v\ra =\la \Pi w,v\ra=\{\Pi w, \Pi v\}\;.
\end{align*}
Furthermore $\{\;,\;\}$ is non-degenerate as $\la\;,\;\ra$ is.
Denote by $\mathcal{Q}$ the closed subscheme of $\mathbb{P}(\Pi \overline{\Lambda})$ consisting of the isotropic lines (with respect to $\{\;,\;\}$). A  subspace $F_0=\la a_1 \Pi e_1 + \ldots + a_n \Pi e_n \ra$ is isotropic if and only if  $\sum  a_i ^2 =0$\,. As $k$ is a field with characteristic $\not=2$ and as $n \geq 3$\,, this implies that $\mathcal{Q}$ is a smooth, geometrically irreducible variety of dimension $n-2$\,.\\
Note that $F_0$ is isotropic with respect to $\{\;,\;\}$ if and only if 
\begin{align*}
F_0 \subseteq (\Pi^{-1}(F_0))^{\bot}\;;
\end{align*} 
with $ ^\bot$ indicating the orthogonal complement with respect to $\la\;,\;\ra $\,. Set 
\begin{align*}
Z_2= \varphi ^{-1}(\mathcal{Q})= \mathcal{Q}\times_{\mathbb{P}(\Pi \overline{\Lambda})} \M_s \;.
\end{align*}
\begin{cl}
$Z_2$ is a locally trivial  $\mathbb{P}^1$-bundle over $\mathcal{Q}$\,.
\end{cl}
\begin{proof}
Denote by $G_0 \subset \Pi \Lambda_{\mathcal{Q}} \subset \Lambda_{\mathcal{Q}}$ the universal $\mathcal{O}_{\mathcal{Q}}$-module corresponding to $\text{id}_{\mathcal{Q}}$\,. It follows that for every $\mathcal{Q}$-scheme $S$ the $S$-valued points of $Z_2$ are all pairs $((G_0)_S,F) \in \mathcal{M}(S)$\,. Set $\mathcal{V}=\Pi^{-1}G_0$\,. Then $\mathcal{V}$ and $\mathcal{V}^{\bot}$ are locally direct summands of  $\Lambda_{\mathcal{Q}}$ of rank $n+1$ and $n-1$\,, respectively. \\
For every pair  $((G_0)_S,F) \in Z_2(S)$\,, the sheaf $F$ satisfies
\begin{align}
\label{z2}
\mathcal{V}^{\bot} _S \subset F \subset \mathcal{V}_S
\end{align}  
On the other hand, if $F$ satisfies (\ref{z2}) we have $\Pi F \subseteq (G_0)_S$ and  $(G_0)_S \subseteq F$ (recall that  $G_0 \subseteq \mathcal{V}^{\bot}$). Furthermore, $F$ is totally isotropic. Indeed, $F/\mathcal{V}^{\bot}_S$ is a submodule of rank $1$ of $\mathcal{V}_S/\mathcal{V}^{\bot} _S$ hence it is isotropic with respect to the induced symplectic form $\la\;,\;\ra$ on   $\mathcal{V}_S/\mathcal{V}^{\bot} _S \times \mathcal{V}_S/\mathcal{V}^{\bot} _S$\,. Set $\mathcal{E}=\mathcal{V}/\mathcal{V}^{\bot}$\,. This is a locally free sheaf of rank $2$ on $\mathcal{Q}$\,. It follows that 
\begin{align*}
Z_2(S)=\{F \subset g^*\mathcal{E}\mid F\text{ is a locally direct summand of rank }1\}\, ,
\end{align*}
where $g$ is the morphism $S \to \mathcal{Q}$\,. Hence $Z_2 \simeq \mathbb{P}(\mathcal{E})$\,.
\end{proof}
\end{st}
\begin{st}
Obviously, $Z_1$ is not contained in $Z_2$ or vice versa. 
\begin{cl}
$\M_s=Z_1 \cup Z_2$\,.
\end{cl}
\begin{proof}
We show that
\begin{align*}
\M_s\setminus Z_2 \rightarrow \mathbb{P}(\Pi \overline{\Lambda})\setminus
 \mathcal{Q}
\end{align*}
is an isomorphism. Choose a point $z \in \mathbb{P}(\Pi \overline{\Lambda})\setminus \mathcal{Q}$\,. Hence $z$ corresponds to a  $\kappa(z)$-subvectorspace $F_0$ of $\Pi \overline{\Lambda}_{\kappa(z)}=\Pi \overline{\Lambda} \otimes \kappa(z)$ of dimension one. Let $S$ be a scheme over $\Spec (\kappa(z))$\,. For every pair $((F_0)_S, F) \in \varphi^{-1}(z)(S)$  we have the relations
\begin{align}
\label{qq}
(F_0)_S \subseteq F \supseteq  (\Pi^{-1}(F_0))_S.
\end{align}
Now $F_0 \not \subset \Pi^{-1}(F_0)$ (recall that $z$ is not contained in $\mathcal{Q}$). We conclude that $(F_0 \oplus \Pi^{-1}(F_0))_S$ is a submodule of $F$\,. But for every point $s \in S$ we have an inclusion $(F_0 \oplus \Pi^{-1}(F_0))_{\kappa(s)}\subset F_{\kappa(s)}$  of vectorspaces of the same dimension, hence they are equal. By Nakayama's lemma, $(F_0 \oplus \Pi^{-1}(F_0))_S=F$\,. On the other hand, the  relations in (\ref{qq}) are also true if we replace $F$ by $\Pi \overline{\Lambda}_S$\,, hence $F= \Pi \overline{\Lambda}_S$\,. It follows that the fiber in $z$ is isomorphic to $\Spec (\kappa(z))$\,. The same argument as in the proof of proposition \ref{red} shows that  the map $\M_s\setminus Z_2 \rightarrow \mathbb{P}(\Pi \overline{\Lambda})\setminus
 \mathcal{Q}$ is an isomorphism.
\end{proof}
\end{st}
\begin{st}
We have shown that  $Z_1$ and $Z_2$ are the irreducible components of $\mathcal{M}_s$\,. The scheme $Z_1 \cap Z_2=Z_1 \times_{\Mls} Z_2$ is smooth of dimension $n-2$ as for every $\Spec (k)$-scheme $S$ we have
\begin{align*}
(F_0,F) \in Z_1(S)\cap Z_2(S) \iff
F=\Pi \overline{\Lambda}_S \text{ and } F_0 \in \mathcal{Q}(S)\;.
\end{align*}
Thus $Z_1 \cap Z_2 \simeq \mathcal{Q}$\,.\\
It follows from  proposition \ref{red}  that $\Ms$ is reduced outside the closed subscheme $Z_1$\,, as $\Mls$ is even smooth outside the point $y$ (see theorem \ref{pappas}). To show that $\Ms$ is reduced even along $Z_1$\,,  we take an open affine covering of $\mathbb{P}(\Lambda)\times Gr_n(\Lambda)$ and write down the equations that define $\mathcal{M}$ in a neighbourhood of $Z_1$\,.\\
We represent $F$ by giving $n$ generating vectors $v_1,\ldots,v_n$ (with respect to the basis $e_1,\ldots,e_n,\Pi e_1,\ldots,\Pi e_n$ of $\Lambda$). Hence we may view $F$ as a $2n \times n$-matrix. The condition that $rk(F)=n$ is expressed by imposing that a certain $n$-minor of the matrix is invertible. As we are only interested in points $F$ with $\overline{F}=\Pi \overline{\Lambda}$\,, we can restrict to the case that 
\begin{align*}
F=&\begin{pmatrix}
a_{11}&\hdots &a_{1n}\\
\vdots&\ddots  &\vdots \\
a_{n1}&\hdots &a_{nn}\\
1&\hdots &0 \\
 \vdots&\ddots  &\vdots \\
0&\hdots &1
\end{pmatrix}\;.
\end{align*}
In the same way we identify $F_0$ with 
\begin{align*}
\begin{pmatrix}
b_1\\
\vdots\\
b_{2n}
\end{pmatrix}
\end{align*}
with the additional condition that $b_k=1$ for some $k \in \{1,\ldots,2n \}$\,. Set $A=(a_{ij})$\,. \\
The condition that $F$ is isotropic translates to 
\begin{align}
\label{e1}
A=A^{t}\;.
\end{align}
The condition $F_0 \subset F$ translates to 
\begin{align}
\label{e*}
\exists \lambda_i\,,i=1 \ldots n : b= \sum_{l=1} ^{n} \lambda_l v_l\;.
\end{align}
But this implies that $\lambda_i=b_{n+i}$ for $i \in \{1,\ldots,n\}$\,. Set
\begin{align*}
\lambda=
\begin{pmatrix}
b_{n+1}\\
\vdots\\
b_{2n}
\end{pmatrix}\;.
\end{align*}
It follows that  equation (\ref{e*}) is equivalent to
\begin{align}
\label{e2}
A\lambda=
\begin{pmatrix}
b_1\\
\vdots\\
b_{n}
\end{pmatrix}\;.
\end{align}
The condition  $(\Pi - \pi)F_0 =(0)$ is equivalent to 
\begin{align}
\label{e3}
A\lambda =\pi \lambda\;.
\end{align}
Note that (\ref{e3}) implies that $b_i\not=1$ for $i=1,\ldots,n$\,. Hence it suffices to consider the case $k \in \{n+1,\ldots, 2n\}$\,. Next we have $(\Pi + \pi)F \subset F_0$\,. This is equivalent to 
\begin{align}
\label{e4}
\exists \gamma=\begin{pmatrix}
\gamma_1\\
\vdots\\
\gamma_n
\end{pmatrix}: A +\pi Id=\lambda \gamma^{t}\;.
\end{align}
The equations (\ref{e2}) and (\ref{e3}) imply that $F$ is $\Pi$-stable. We may replace $A$ by $\lambda \gamma^{\bot} - \pi Id$\,. Hence (\ref{e1}) -- (\ref{e4}) is equivalent to
\begin{align}
\label{e5}
\lambda \gamma^{t} &= \gamma \lambda^{t}\;;\\
\label{e6}
\lambda \gamma^{t} \lambda &= 2\pi \lambda \;;\\
\label{e7}
\lambda_k&=1\;.
\end{align}
Using (\ref{e7}), equation (\ref{e5}) is equivalent to
\begin{align*}
\gamma_i= \gamma_k \lambda_i \,, i=1,\ldots,n\;.
\end{align*}
Hence for $i\not=k, \gamma_i$ is determined by $\lambda_i$ and $\gamma_k$\,. It follows that equation (\ref{e6}) translates to 
\begin{align*}
(\sum_{i=1} ^n \lambda_i ^2)\gamma_k=2\pi\;.
\end{align*}
Set $y=\gamma_k$ and $x_i=\lambda_i$\,. We have shown that the open subschemes  \begin{align}\label{gleichung}
U_k &= \Spec (R[x_1,\ldots,x_n,y]/(y (\sum x_i ^2) - 2\pi, x_k -1))
\end{align}
of $\M$  cover $Z_1$\,. We have 
\begin{align*}
Z_1\cap U_k & = \Spec (k[x_1,\ldots,x_n]/(x_k -1))\;;\\
Z_2 \cap U_k & = \Spec (k[x_1,\ldots,x_n,y]/(x_k -1, \sum x_i ^2))\;; \\
(Z_1 \cap Z_2) \cap U_k & = \Spec (k[x_1,\ldots,x_n]/(x_k -1, \sum x_i ^2))\;.
\end{align*}
We conclude that $\mathcal{M}_s$ is reduced.
\end{st}
\begin{st}
The morphism $\mathcal{M} \rightarrow \Spec (R)$ is flat.
\begin{proof}
If $x$ is a point of $\M$ which does not lie on $Z_1$\,, then the morphism  $\mathcal{M} \rightarrow \Spec (R)$ is flat in $x$\,, as $\M \setminus Z_1 \rightarrow \Ml\setminus y$ is an isomorphism (see proposition \ref{red}) and the latter scheme is even smooth over $\Spec(R)$ (see theorem \ref{pappas}). If $x \in Z_1$\,,  then there is a number $k$ such that  $x \in U_k$ as in (\ref{gleichung}). But it is obvious that $U_k$ is flat over $\Spec(R)$\,.\end{proof}
\end{st}
\begin{st}
It remains to show that $\mathcal{M}$ is  regular of dimension $n$\,. Since the generic fiber $\M_\eta$ of $\M$ is smooth, the local ring $\mathcal{O}_{\mathcal{M},x}$ is regular of dimension $n-1$ if $x \in \M_{\eta}$ is a closed point (recall that $\dim \M_{\eta}=n-1$). As $\mathcal{M} \rightarrow \Spec (R)$ is flat, it follows that for all closed points $x \in \mathcal{M}$ lying on the special fiber we have $\dim \mathcal{O}_{\mathcal{M},x}=n$ and that $\mathcal{M}$ is regular in $x$ if $x$ is not contained in  $Z_1 \cap Z_2$\,, because $\M$ is smooth outside $Z_1 \cap Z_2$\,. Hence we only have to show that $ \mathcal{O}_{\mathcal{M},x}$ is regular for closed points  lying on $  Z_1 \cap Z_2$\,. Let $x$ be such a point. There is a number $k$ such that $x \in U_k$\,. We have 
\begin{align*}
\OO_{Z_1 \cap Z_2,x}&=\OO_{\M,x}/(y, \sum x_i ^2,\pi) \\
&=\OO_{\M,x}/(y, \sum x_i ^2)
\end{align*}
(as  $2\pi=0$ implies  $\pi=0$). As $B=\OO_{Z_1 \cap Z_2,x}$ is a regular local ring of dimension $n-2$\,, there are elements $\overline{t}_1,\ldots,\overline{t}_{n-2}$ generating $m_B/{m}_B ^2$\,. We can lift these elements to  $C=\OO_{\M,x}$\,. Hence $t_1,\ldots,t_{n-2},y, \sum x_i ^2$ generate $m_C/{m}_C ^2$\,. On the other hand $\dim \OO_{\M,x}=n$\,, hence $\M$ is regular in $x$\,.
\end{st}
As the dimension of the generic fiber is $n-1$ and for all closed points $x\in \mathcal{M}$ lying on the special fiber $\dim \mathcal{O}_{\mathcal{M},x}=n$\,, the dimension of $\mathcal{M}$ is $n$\,. This completes the proof of theorem \ref{thm:reg}\,.
\end{proof}
\section[Trace of Frobenius]{Trace of Frobenius on the sheaf of nearby cycles}\label{frob}
Denote by $R$ a henselian discrete valuation ring. Let $X$ be a scheme of finite type over $S=\Spec (R)$\,. Let $\overline{\eta}$ be a geometric point over $\eta$ and denote by $\overline{s}$ the corresponding geometric point over $s$\,. Denote the characteristic exponent of $s$ by $p$\,. Consider the diagram
\begin{align*}
\begin{CD}
X_{\overline{\eta}} @>\overline{j}>> X @<{\tau}<< X_s\\
@VVV   @VVV @VVV\\
\overline{\eta} @>>> S @<<< s
\end{CD}\;.
\end{align*}
\begin{df}
Let $l$ be a prime number with $l \not= p$\,. The sheaf
\begin{eqnarray*}
R^i\Psi \mathbb{Q}_l=\tau^* R^i \overline{j}_* \mathbb{Q}_l
\end{eqnarray*}
on $X_s$ is called the $i$-th sheaf of nearby cycles. Denote by $I \subseteq \text{Gal}(\overline{\eta}/\eta)$ the inertia group and by $P \subseteq I$ the biggest subgroup of $I$ which is a pro-$p$-group. The sheaf
\begin{align*}
R^i\Psi_t \mathbb{Q}_l=(R^i\Psi \mathbb{Q}_l)^P
\end{align*}
is called the $i$-th sheaf of tame nearby cycles.
\end{df}
Suppose that the residue field of $s$ is contained in $\mathbb{F}_q$\,. Consider the geometric Frobenius $\frob_q$ over $\mathbb{F}_q$ i.e.\ the map $x\mapsto x^{1/q}$ on $\kappa(\overline{s})$\,. For a  $\mathbb{F}_q$-valued point $x$ of $X_s$\,, we have the alternating semisimple trace

\begin{align*}
\Tr^{\text{ss}}(\frob_q;R\Psi_t(\mathbb{Q}_l)_x ^I)= \sum_i (-1)^i
\Tr^{\text{ss}}(\frob_q;R^i \Psi_t(\mathbb{Q}_l)_x ^I)\;.
\end{align*}

We are going to compute this alternating semisimple trace in the case where $X$ is the local model $\Ml$ introduced in section \ref{local model} and $(r,s)=(n-1,1)$\,. We do this in two steps: First we calculate the trace on $\M$ and then relate it to the trace on $\Ml$\,. We assume in the rest of the section that $k=R/\pi$ is contained in $\mathbb{F}_q$\,.
\begin{thm}
\label{thm:nbc}
Let $X$ be a scheme over $S$ of finite type with smooth generic fiber. Suppose that the special fiber is a reduced divisor with normal crossings. Assume further that $X_{\overline{s}}$ is globally the union of smooth irreducible divisors. Let $x\in X_s$ and let $S_x$ be the set of irreducible components of $X_{\overline{s}}$ passing through $x$\,. Then
\begin{align*}
R^1\Psi_t (\mathbb{Q}_l)_x&=\text{ker}(\bigoplus_{S_x}\mathbb{Q}_l(-1) \stackrel{\Sigma}{\rightarrow} \mathbb{Q}_l(-1))\;; \\
R^i\Psi_t(\mathbb{Q}_l)_x&=\sideset{}{^i}\bigwedge  R^1\Psi_t (\mathbb{Q}_l)_x\;.
\end{align*}
In particular, the inertia group acts trivially on $R^i\Psi_t (\mathbb{Q}_l)_x$\,.
\end{thm}
\begin{proof}
See \cite{SGA7}, exp. I, theorem  3.3\,.
\end{proof}
The operation of ${\frob}_q$ on the $\mathbb{Q}_l$-vectorspace $V:=\mathbb{Q}_l(-1)^n$ is just multiplication with $q$\,,  hence we have $\Tr^{\text{ss}}({\frob}_q,V)=nq$\,. It follows that for $j \in \{1,\ldots,n\}$ we have $\Tr^{\text{ss}}(\frob_q,\bigwedge ^j V)=\binom{n}{j}q^j$\,.\\
As the scheme $\mathcal{M}$ satisfies the conditions of theorem
\ref{thm:nbc} (recall that we assume that $k$ is contained in $\mathbb{F}_q$)\,, we have
\begin{align*}
\Tr^{\text{ss}}(\frob_q,R\Psi_t(\mathbb{Q}_l)_x ^I)= \sum_{j=0} ^t (-1)^j \binom{t}{j}q^j = (1-q)^t
\end{align*}
for every $x \in \mathcal{M}_s$ (with $t=|S_x| -1$). The following theorem connects the alternating semisimple traces of Frobenius on $\mathcal{M}$ resp. $M^{loc}$\,.
\begin{thm}
Assume that $\kappa(s)$ is contained in $\mathbb{F}_q$ and let $X' \stackrel{f}{\rightarrow}X$ be a proper morphism of $S$-schemes of finite type. Suppose that $f$ is an isomorphism on the generic fibers, and denote by $R \Psi_t$ resp. $R{\Psi_t}'$ the sheaves of tame nearby cycles on $X$ resp. $X'$\,. For $x \in X(\mathbb{F} _q)$ we have
\begin{align*}
\Tr^{\text{ss}}(\frob_q;R\Psi_t(\mathbb{Q}_l)_x)= \sum_{\stackrel{x' \in X'(\mathbb{F}_q)}{x'\mapsto x}}\Tr^{\text{ss}}(\frob_q;R{\Psi_t}'(\mathbb{Q}_l)_{x'})
\end{align*} 
\end{thm}
\begin{proof}
Compare \cite{Ulrich}, corollary 4.3\,.
\end{proof}
In section \ref{resolution} we saw that the morphism 
\begin{align*}
\M(\mathbb{F}_q) \stackrel{\tau(\mathbb{F}_q)}{\longrightarrow} \Ml(\mathbb{F}_q)
\end{align*}
is bijective outside the point  $\Pi\Lambda_{\mathbb{F}_q}=\Pi\Lambda \otimes \mathbb{F}_q$ and that for every point $x \not= \Pi\Lambda_{\mathbb{F}_q}$ the preimage lies on only  one irreducible component (namely  $Z_2$). In this case we have
\begin{align*}
\Tr^{\text{ss}}(\frob_q;R\Psi_t(\mathbb{Q}_l)_x)= 1\;.
\end{align*} 
Now consider the $\Spec (\mathbb{F}_q)$-valued point $y= \Pi\Lambda_{\mathbb{F}_q}$ of $\Ml$\,. As explained in section \ref{resolution} the preimage of $y$ under the map $\M(\mathbb{F}_q) \rightarrow \Ml(\mathbb{F}_q)$ is $\mathbb{P}(\Pi \overline{\Lambda})(\mathbb{F}_q)$\,. By definition, every $\mathbb{F}_q$-rational point point $F_0=\spa \{a_1\Pi e_1 + \ldots +  a_n \Pi e_n\}$ of $\tau^{-1}(y)$ lies on the irreducible component $Z_1$ and it lies on $Z_2$ if and only if $\sum a_i ^2 =0$. Hence the number of $\mathbb{F}_q$-rational points of $\tau^{-1}(y)$  lying on both irreducible components equals
\begin{align*}
z= |\mathcal{Q}(\mathbb{F}_q)| =\frac{|\{(a_1,\ldots,a_n)\in \mathbb{F}_q ^n |\sum a_i ^2 =0\}|-1}{q-1} \;.
\end{align*}
Hence for the point $y$  the alternating semisimple trace of Frobenius equals  
\begin{align*}
\Tr^{\text{ss}}(\frob_q;R\Psi_t(\mathbb{Q}_l)_y)&=z \cdot (1-q) +(|\mathbb{P}(\Pi \overline\Lambda)(\mathbb{F}_q)| -z)\cdot 1\\
&=z \cdot (1-q) + \frac{q^n -1}{q-1}- z\\
&=\frac{q^n -1}{q -1} -z \cdot q \;.
\end{align*}
There is an explicit formula for $z$ (see e.g. \cite{Kitaoka}, lemma 1.3.1):
\begin{align*}
z=\frac{q^{n-1} - \chi q^{\frac{n}{2}-1} + \chi q^{\frac{n}{2}}-1}{q-1}\;;
\end{align*}
with $\chi=0$ if $n$ is odd and $\chi=1$ or $-1$ if $n$ is even, depending on whether $\mathbb{F}_q ^n$ is hyperbolic with respect to $v\mapsto \sum v_i ^2$ or not. Combining these two results we obtain the following theorem.
\begin{thm}
Let $x$ be a $\mathbb{F}_q$-rational point of $\Ml$ and denote by $y$ the point corresponding to $\Pi \overline{\Lambda}$\,. If $x \not=y$ we have
\begin{align*}
\Tr^{\text{ss}}(\frob_q;R\Psi_t(\mathbb{Q}_l)_x)= 1\;.
\end{align*}
For $x=y$ we have 
\begin{align*}
\Tr^{\text{ss}}(\frob_q;R\Psi_t(\mathbb{Q}_l)_y)=
\begin{cases}
1 & \text{if }n \text{ is odd}\\
1 - q^{\frac{n}{2}} & \text{if }n \text{ is even and } \mathbb{F}_q ^n \text{ is hyperbolic}\\
1 + q^{\frac{n}{2}} &  \text{if }n \text{ is even and } \mathbb{F}_q ^n \text{ is not hyperbolic}\;.
\end{cases}
\end{align*}
\end{thm}

\end{document}